\documentclass[a4paper,12pt]{amsart}

\usepackage{amsmath}
\usepackage{amsthm}
\usepackage{amssymb}
\usepackage{amsfonts}
\usepackage{algorithm}
\usepackage{algpseudocode}
\usepackage{bbm}
\usepackage{dsfont}
\usepackage{graphicx}
\usepackage{pstricks,pstricks-add,pst-math,pst-xkey}
\usepackage{hyperref}

\newtheorem{theo}{Theorem}
\newenvironment{theorem}{\vspace{4mm}\begin{theo}}{\end{theo}}
\newtheorem{prop}[theo]{Proposition}
\newenvironment{proposition}{\vspace{4mm}\begin{prop}}{\end{prop}}
\newtheorem{lem}[theo]{Lemma}

\newtheorem{coro}[theo]{Corollary}

\newtheorem{rem}[theo]{Remark}
\newenvironment{remark}{\vspace{4mm}\begin{rem}\rm}{\end{rem}}
\newtheorem{exa}[theo]{Example}

\newtheoremstyle{citing}{}{}{\itshape}{}{\bfseries}{.}%
 { }{\thmnote{#3}}
\theoremstyle{citing}
\newtheorem{cit}{}

\newcommand{\Z}{\mathbb{Z}}
\newcommand{\N}{\mathbb{N}}

\newcommand{\ind}{\scalebox{1.3}{\raisebox{-1.5pt}{$\mathds{1}$}}}
\newcommand{\tvd}[1]{\Vert #1 \Vert_\text{TV}}	
\newcommand{\abs}[1]{\left\vert #1 \right\vert}	

\newcommand{\E}{\mathbb{E}}

\renewcommand{\O}{\Omega}

\begin{document}

\title[Exact Sampling for the Ising model]
	{Exact Sampling for the Ising model at all temperatures}
\author[Mario Ullrich]{Mario Ullrich\\
	Friedrich Schiller University Jena \vspace*{5mm}\\
	{\today}}
\thanks{The author was supported by the DFG GK 1523.}
\address{Friedrich Schiller University Jena\newline
		\textit{homepage:} \url{http://users.minet.uni-jena.de/~ullrich/}}
\email{mario.ullrich@uni-jena.de}




\begin{abstract}
The Ising model is often referred to as the most studied model 
of statistical physics. It describes the behavior of ferromagnetic 
material at different temperatures. 
It is an interesting model also for mathematicians, 
because although the Boltzmann distribution is continuous in the 
temperature parameter, the behavior of the usual single-spin 
dynamics to sample from this measure varies extremely.
Namely, there is a critical temperature where we get 
rapid mixing above and slow mixing below this value. 
Here, we 
give a survey of the known results on mixing time of 
Glauber dynamics for the Ising model on the square lattice and 
present a technique that makes exact sampling of the Ising 
model at all temperatures possible 
in polynomial time. At high temperatures this is well-known 
and 
although this seems to be known also in the low temperature case 
since Kramer and Waniers paper \cite{KW} from the 1950s, 
we did not found any reference that 
describes exact sampling for the Ising model at low temperatures.
\end{abstract}

\maketitle

\thispagestyle{empty}

\section{Introduction}

In this article we summarize the known results about the mixing time 
of the heat bath dynamics for the Ising model and combine them 
with some graph theoretic results to an algorithm to sample 
exactly from the Ising model in polynomial time. 
By time (or running time) we always mean the number of steps of 
the underlying Markov chain.
The algorithm that will be analyzed (Algorithm 2, given 
in Section \ref{sec-efficient}) is at high temperatures simply 
the Coupling from the past algorithm (see Propp and Wilson \cite{PW}).
At low temperatures we have to produce a sample at the dual graph, 
but this can be traced back to sampling on the initial graph with 
constant boundary condition.

The main theorem of this article is stated as follows.

\begin{cit}[Theorem~\ref{th-main}]
Let $G_L$ be the square lattice with $N=L^2$ vertices. Then, 
Algorithm 2 outputs an exactly distributed Ising 
configuration with respect to $\pi_\beta^{G_L}$ in expected time 
smaller than
\begin{itemize}
	\item \quad $c_\beta\, N \,(\log N)^2$ 
		\qquad for $\beta\neq\beta_c=\log(1+\sqrt{2})$ and some $c_\beta>0$
	\vspace{2mm}
	\item \quad $16\, N^C \log N$ 
		\qquad
		for $\beta=\beta_c$, where $C$ is given in \eqref{eq-critical-C}.
\end{itemize}
\end{cit}

As a consequence we get that one can estimate the expectation of 
arbitrary functions with respect to the Boltzmann distribution 
in polynomial time. Namely, if we use the simple Monte Carlo method 
to approximate the expectation of a function $f$ on the Ising model, 
we need $\epsilon^{-2}\Vert f\Vert^2_{2}$ exact samples from 
$\pi_\beta$ (i.e. Algorithm~2) to reach a mean square error of at 
most $\epsilon$. 
Therefore, if we denote the bounds from Theorem~\ref{th-main} by 
$T_\beta$, we need on average 
$T_\beta\,\epsilon^{-2}\Vert f-\E_{\pi_\beta}f\Vert^2_{2}$ 
steps of the Markov chain that will be defined in Section 
\ref{sec-ising} \\
The first polynomial-time algorithm (FPRAS) was shown by Jerrum and 
Sinclair \cite{JS}. There they present an algorithm to approximate the 
partition function $Z_\beta$ and, as a consequence, approximate 
expectations of functions that are given in terms of the partition 
function in polynomial time at all temperatures $\beta$.


\vspace{3mm}
\section{The Ising model}\label{sec-ising}

In this section we introduce the two-dimensional Ising model. \\
Let $G=(V,E)$ be a graph with finite vertex set $V\subset\Z^2$ and 
edge set $E=\left\{\{u,v\}\in\binom{V}{2}:\, \abs{u-v}=1 \right\}$, 
where $\binom{V}{2}$ is the set of all subsets of $V$ with 2 
elements. From now, $N:=\abs{V}$.
We are interested in the square lattice, i.e. $V=\{1,\dots,L\}^2$ 
for some $L=\sqrt{N}\in\N$, because it is the most widely used 
case. We denote the induced graph by $G_L$.\\
The \emph{Ising model} on $G_L$ is now defined as the set of possible 
configurations $\O_{\rm IS}=\{-1,1\}^V$, where $\sigma\in\O_{\rm IS}$ 
is an assignment of -1 or 1 to each vertex in $V$, 
together with the probability measure
\vspace{1mm}
\[
\pi_\beta(\sigma) \;:=\; \pi^{G_L}_\beta(\sigma) \;=\; \frac1{Z_\beta}\,
	\exp\left\{\beta\,\sum_{u,v:\, u\leftrightarrow v}
	\ind\bigl(\sigma(u)=\sigma(v)\bigr)\right\},
\]
where $u\leftrightarrow v$ means $u$ and $v$ are neighbors in $G_L$, 
$Z$ is the 
normalization constant and $\beta\ge0$ is the called the inverse 
temperature. This measure is called the Boltzmann (or Gibbs) 
distribution with free boundary condition.\\
Additionally we need the notion of \emph{boundary conditions}, but we 
restrict ourself here to the ``all plus'' and ``all minus'' case.\\
Let $V^c=\Z^2\setminus V$.
Then we denote the lattice $G_L$ together with the probability 
measure
\[
\pi_\beta^{\pm}(\sigma) \;:=\; \pi^{G_L,\pm}_\beta(\sigma) 
	\;=\; \frac{1}{\widetilde Z_\beta}\; \pi^{G_L}_\beta(\sigma)\cdot
	\exp\left\{\beta\,
	\sum_{\substack{v\in V, \,u\in V^c:\\ u\leftrightarrow v} } 
	\ind\Bigl(\sigma(v)=\pm1\Bigr)\right\}
\]
by the Ising model with plus/minus boundary condition, respectively. 
One can imagine that this corresponds to the Ising model on $G_L$ with 
a strip of fixed spins around, so every vertex in $G_L$ has the same 
number of neighbors.\\
In 1944 Onsager \cite{Onsager} proved that there is a phase transition 
at $\beta=\beta_c:=\ln(1+\sqrt{2})$ in the case where $V=\Z^2$ and 
we will see that this value is also important for finite 
lattices. Namely, the dynamics that will be defined below is 
rapidly mixing if and only if $\beta\le\beta_c$.\\
We will use the so called \emph{heat bath dynamics}. 
These dynamics define a irreducible, aperiodic and reversible 
Markov chain $X^\beta=(X_i^\beta)_{i\in\N}$ with 
stationary distribution $\pi_\beta$ by the transition matrix
\vspace{1mm}
\[
P(\sigma,\sigma^{v,\xi}) \;=\; \frac{1}{N}\;
	\left(1+\frac{\pi_\beta(\sigma)}{\pi_\beta(\sigma^{v,\xi})}\right)^{-1},
	\qquad \sigma\in\O_{\rm IS},\; v\in V,
\vspace{1mm}
\]
where $\sigma^{v,\xi}$ with $\xi\in\{-1,1\}$ is defined by 
$\sigma^{v,\xi}(v)=\xi$ and $\sigma^{v,\xi}(u)=\sigma(u)$, $u\neq v$. 
The interpretation of this algorithm is very simple. In each step 
choose a random $v\in V$ and assign a new value to $v$ according 
to $\pi_\beta$ conditioned on all the neighbors of $v$.\\
Note that the results of this article hold in general for all 
Glauber dynamics as defined in \cite{Glauber} that admit a monotone 
coupling (see Section \ref{sec-sampling}). For a general introduction 
to Markov chains see e.g. \cite{LPW}, or \cite{M} in the context of 
spin systems.

In the sequel we want to estimate how fast such a Markov chain 
converges to its stationary distribution. Therefore we first introduce 
the \emph{total variation distance} to measure the distance between 
two probability measures $\nu$ and $\pi$, which is defined by
\[
\tvd{\nu-\pi} \;=\; \frac12\,\sum_{\sigma\in\O_{\rm IS}}\,
	\abs{\nu(\sigma)-\pi(\sigma)}.
\]
Now we can define the \emph{mixing time} of the Markov chain with 
transition matrix $P$ and stationary distribution $\pi_\beta$ by 
\vspace{1mm}
\[
\tau_\beta \;=\; \min\left\{n: \max_{\sigma\in\O_{\rm IS}}
\tvd{P^n(\sigma,\cdot)-\pi_\beta(\cdot)}\,\le\,\frac1{2\rm e}\right\}.
\vspace{1mm}
\]
This is the expected time the Markov chain needs to get close to 
its stationary distribution. In fact, one can bound the spectral gap 
of the transition matrix $P$ in either direction in terms of the 
mixing time, see e.g. \cite[Th. 12.3 \& 12.4]{LPW}, so one can bound 
the error of a MCMC algorithm to integrate functions over $\O_{\rm IS}$, 
as one can read in \cite{Rud}. 
Furthermore, if the Markov chain is rapidly mixing 
(i.e. the mixing time is at most polylogarithmic in the 
size of the state space $\O_{\rm IS}$) we get that the problem 
of integration (with an unnormalized density) on the Ising model
is \emph{tractable}, see also \cite{NW2}. Unfortunately, there is no 
Markov chain that is proven to be rapidly mixing at all temperatures.\\
However, in this article we are interested in sampling exactly from 
the stationary distribution, but first we present the known 
mixing time results for the Glauber dynamics for the Ising model.
For proofs or further details we refer to the particular articles 
or the survey of Martinelli \cite{M}.
Of course, we can only give a small selection of references, 
because there are many papers leading to the results given below. 

\vspace{2mm}
\begin{theorem}{\cite{MO1}}\label{th-mix-high}
Let $\beta<\beta_c$. Then there exists a constant $c_\beta>0$ such 
that the mixing time of the Glauber dynamics 
for the Ising model with arbitrary boundary condition on 
$G_L$ satisfies
\[
\tau_\beta \;\le\; c_\beta\;N\log N.
\]
\end{theorem}

\begin{theorem}{\cite{CGMS}}\label{th-mix-low}
Let $\beta>\beta_c$. Then there exists a constant $c_\beta>0$ such 
that the mixing time of the Glauber dynamics for the Ising model 
on $G_L$ satisfies
\[
\tau_\beta \;\ge\; e^{c_\beta N}.
\]
\end{theorem}
\vspace{2mm}

The results above can be obtained by the observation that 
some spatial mixing property of the measure $\pi_\beta$ is 
equivalent to the mixing in time of the Glauber dynamics. 
For details for this interesting fact, see \cite{DSVW}.\\
The constant $c_\beta$ of Theorem \ref{th-mix-high} is 
widely believed to be of order 
$\frac1{\beta-\beta_c}$. To determine the mixing time in the 
case $\beta=\beta_c$ was a challenging problem for a long time. 
It was solved by Lubetzky and Sly in their recent paper \cite{LS}.

\begin{theorem}{\cite{LS}}\label{th-mix-critical}
There exists a constant $C>0$ such that the mixing time of the 
Glauber dynamics for the Ising model 
on $G_L$ at the critical temperature satisfies
\[
\tau_\beta \;\le\; 4\,N^C.
\]
\end{theorem}

\begin{remark}
We give here only a brief description of the constant $C$, which 
can be given explicitly. For more details see \cite[p.19]{LS}.\\
However, numerical experiments on the ``true'' exponent suggest that 
$C\approx3.08$ (see e.g. \cite{WHS}, \cite{NB} and note the explanation 
below). \\
The constant $C$ in Theorem \ref{th-mix-critical} is given by 
\begin{equation}\label{eq-critical-C}
C\;=\;2+\log_{3/2}\left(\frac{2}{1-p^+}\right).
\end{equation}
Here, $p^+$ is the limiting vertical crossing probability in 
the random cluster model on a fully-wired rectangle, where
the width of the lattice is 3 times its height.
The $C$, as given here, differs from the one given in \cite{LS} 
by eliminating a factor of 2 in front of the $\log$ term and by 
the additional 2. The reason is that we state their result in 
terms of $N$ and not in the side-length $L$ of the lattice
(therefore without factor 2) and that we are interested in the 
discrete time single-spin algorithms. Therefore we get an additional 
factor $N$ in their spectral gap result (\cite[Th. 1]{LS}) and a 
factor $N$ by (see e.g. \cite{LPW})
\[
\tau_\beta\;\le\;\log\left(\frac{e}{\min_{\sigma}\pi_{\beta_c}(\sigma)}\right)
	\,\text{\rm\bf gap}(X^\beta)^{-1}
\;\le\; 4\,N\,\text{\rm\bf gap}(X^\beta)^{-1},
\]
because $\min_{\sigma}\pi_{\beta_c}(\sigma)\ge \exp(-3N)$.
\end{remark}

The results of this section show that the Glauber dynamics is 
rapidly mixing 
for $\beta\le\beta_c$, but very slowly mixing for larger $\beta$. 
In Section \ref{sec-rc} we will see how to avoid this problem.

\section{Exact sampling}\label{sec-sampling}

In this section we briefly describe the so called 
\emph{Coupling from the past algorithm} (CFTP) to sample exactly 
from the stationary distribution of a Markov chain.\\
This algorithm works under weak assumptions on the Markov 
chain for every finite state space and every distribution, but to 
guarantee that the algorithm is efficient we need some monotonicity 
property of the model and that the chain is rapidly mixing. 
For a detailed description of CFTP and the proof of 
correctness see \cite{PW}.\\
We restrict ourself to the heat bath dynamics for the Ising model. 
First note that the heat bath dynamics, as defined above, admits 
a monotone coupling, that is, given two realizations of the heat bath 
chain $X=(X_t)_{t\in\N}$ and $Y=(Y_t)_{t\in\N}$, there exists a 
coupling $(X,Y)$ (i.e. using the same random numbers) such that
\[
X_t \;\le\; Y_t \;\;\Longrightarrow\;\; X_{t+1} \;\le\; Y_{t+1} 
\qquad \text{ for all } t\in\N,
\]
where $\le$ means smaller or equal at each vertex.\\
Additionally we know that $-\bf{1}\le\sigma\le\bf{1}$ for all 
$\sigma\in\O_{\rm IS}$, where ${-\bf{1}}=(-1)^V$ and ${\bf{1}}=(1)^V$.
Therefore if we set $X_0=-{\bf1}$ and $Y_0={\bf1}$ we know that 
$X_0\le\sigma\le Y_0$ for all $\sigma$ and so 
$X_t\le Z_t\le Y_t$ for the realization $Z=(Z_t)_{t\in\N}$ 
with $Z_0=\sigma$. Since this holds for all $\sigma$, one can choose 
$Z_0\sim\pi_\beta$ and we get that whenever $X_t$ and $Y_t$ coalesce, 
they also coalesce with $Z_t$ which has the right distribution.\\
After we presented the idea of the algorithm, we state the algorithm 
in detail. Note that the algorithm is called Coupling from the past, 
because we run the chains from the past to the present.
The algorithm $\text{CFTP}(G,\beta)$ to sample from the distribution 
$\pi^G_\beta$ works as described in Algorithm 1.

\begin{algorithm}
\caption{\quad Coupling from the past}
\begin{algorithmic}[1]
\Statex\Call{\bf Input:}{} The graph $G=(V,E)$ and the value of $\beta$
\Statex\Call{\bf Output:}{} An Ising configuration $\sigma\sim\pi_\beta$
\vspace{2mm}
\Procedure{CFTP}{$G,\beta$}
\vspace{2mm}
\State Set $t = 0$
\State Set $X_0=-{\bf1}$ and $Y_0={\bf1}$
\vspace{2mm}
\While{$X_0 \neq Y_0$}
\State $t = t+1$
\vspace{1mm}
\State Generate random numbers $U_{-2^t+1},\dots,U_{-2^{t-1}}$ that are 
\Statex	\qquad\qquad sufficient to run the Markov chain. 
\Statex	\qquad\qquad\quad (e.g. $U_i\sim \text{Uniform }\{V\times[0,1]\}$)
\vspace{1mm}
\State Set $X_{-2^t+1}={\bf0}$ and $Y_{-2^t+1}={\bf1}$ and run the chains 
					until 
\Statex \qquad\qquad time 0 by using only the random numbers 
								$U_{-2^t+1},\dots,U_{-1}$
\EndWhile
\vspace{1mm}
\State \textbf{return} $\sigma = X_0$
\vspace{1mm}
\EndProcedure
\end{algorithmic}
\end{algorithm}

We denote the algorithm by $\text{CFTP}^\pm(G,\beta)$ if we sample 
with respect to $\pi^\pm_\beta$, i.e. with plus/minus boundary 
condition.

See \cite{H} for examples that show that it is necessary to go 
from the past in the future and that we have to reuse the random 
numbers.\\
Now we state the connection between the expected running time of 
the CFTP algorithm and the mixing time of the Markov chain.

\vspace*{5mm}
\begin{proposition}{\cite{PW}}\label{prop-cftp}
Let $T_\beta$ be the expected running time of CFTP$(G,\beta)$ from 
Algorithm 1 with $G=(V,E)$ and $\abs{V}=N$. Then 
\[
T_\beta \;\le\; 4\, \tau_\beta \,\log N,
\]
where $\tau_\beta$ is the mixing time of the underlying Markov chain.
\end{proposition}

We see that exact sampling from the Boltzmann distribution is efficient 
whenever the Markov chain is rapidly mixing.
By the results of Section \ref{sec-ising} we know that this is the case 
for $\beta\le\beta_c$. In the case $\beta>\beta_c$ we need a different 
technique to generate exact samples. Therefore we need essentially 
the so called random cluster model, as we will see in the next section.

\section{The random cluster model}\label{sec-rc}

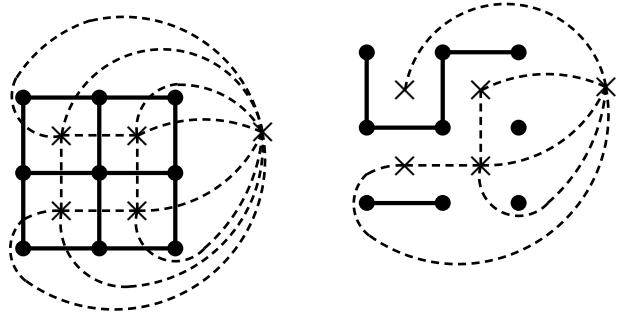
\begin{figure}[ht]
\scalebox{1}{\psset{xunit=1.0cm,yunit=1.0cm,algebraic=true,dotstyle=*,dotsize=3pt 0,linewidth=0.8pt,arrowsize=3pt 2,arrowinset=0.25}
\begin{pspicture*}(-0.43,-0.06)(6.02,4.43)
\psline[linewidth=1.6pt](1,3)(2,3)
\psline[linewidth=1.6pt](2,3)(3,3)
\psline[linewidth=1.6pt](3,3)(3,2)
\psline[linewidth=1.6pt](3,2)(3,1)
\psline[linewidth=1.6pt](3,1)(2,1)
\psline[linewidth=1.6pt](2,1)(2,2)
\psline[linewidth=1.6pt](2,2)(3,2)
\psline[linewidth=1.6pt](2,2)(2,3)
\psline[linewidth=1.6pt](1,3)(1,2)
\psline[linewidth=1.6pt](1,2)(2,2)
\psline[linewidth=1.6pt](2,1)(1,1)
\psline[linewidth=1.6pt](1,1)(1,2)
\psline[linewidth=1.0pt,linestyle=dashed,dash=3pt 2pt](1.5,1.5)(1.5,2.5)
\psline[linewidth=1.0pt,linestyle=dashed,dash=3pt 2pt](1.5,2.5)(2.5,2.5)
\psline[linewidth=1.0pt,linestyle=dashed,dash=3pt 2pt](2.5,2.5)(2.5,1.5)
\psline[linewidth=1.0pt,linestyle=dashed,dash=3pt 2pt](2.5,1.5)(1.5,1.5)
\parametricplot[linewidth=1.0pt,linestyle=dashed,dash=3pt 2pt]{0.1852913062710569}{2.989280090214289}{1*1.34*cos(t)+0*1.34*sin(t)+2.83|0*1.34*cos(t)+1*1.34*sin(t)+2.3}
\parametricplot[linewidth=1.0pt,linestyle=dashed,dash=3pt 2pt]{1.1350083349985625}{2.0596018113318753}{1*1.84*cos(t)+0*1.84*sin(t)+3.37|0*1.84*cos(t)+1*1.84*sin(t)+0.87}
\parametricplot[linewidth=1.0pt,linestyle=dashed,dash=3pt 2pt]{4.6717617082515}{5.883637174045906}{1*1.71*cos(t)+0*1.71*sin(t)+2.57|0*1.71*cos(t)+1*1.71*sin(t)+3.21}
\parametricplot[linewidth=1.0pt,linestyle=dashed,dash=3pt 2pt]{1.5028552028831879}{2.100187227182034}{1*0.86*cos(t)+0*0.86*sin(t)+1.44|0*0.86*cos(t)+1*0.86*sin(t)+0.64}
\parametricplot[linewidth=1.0pt,linestyle=dashed,dash=3pt 2pt]{2.311045279311494}{4.0371577644321395}{1*0.54*cos(t)+0*0.54*sin(t)+1.37|0*0.54*cos(t)+1*0.54*sin(t)+0.99}
\parametricplot[linewidth=1.0pt,linestyle=dashed,dash=3pt 2pt]{-2.1993684311435224}{0.19067988888988038}{1*1.98*cos(t)+0*1.98*sin(t)+2.2|0*1.98*cos(t)+1*1.98*sin(t)+2.17}
\parametricplot[linewidth=1.0pt,linestyle=dashed,dash=3pt 2pt]{2.8453007051607155}{5.540388001593068}{1*0.52*cos(t)+0*0.52*sin(t)+3|0*0.52*cos(t)+1*0.52*sin(t)+1.35}
\parametricplot[linewidth=1.0pt,linestyle=dashed,dash=3pt 2pt]{5.428195851138397}{6.21778798060977}{1*2.24*cos(t)+0*2.24*sin(t)+1.91|0*2.24*cos(t)+1*2.24*sin(t)+2.69}
\parametricplot[linewidth=1.0pt,linestyle=dashed,dash=3pt 2pt]{1.651943081227404}{3.288962910742752}{1*0.59*cos(t)+0*0.59*sin(t)+3.08|0*0.59*cos(t)+1*0.59*sin(t)+2.59}
\parametricplot[linewidth=1.0pt,linestyle=dashed,dash=3pt 2pt]{0.523108668073518}{1.593136560077821}{1*1.25*cos(t)+0*1.25*sin(t)+3.06|0*1.25*cos(t)+1*1.25*sin(t)+1.92}
\parametricplot[linewidth=1.0pt,linestyle=dashed,dash=3pt 2pt]{2.956350314060717}{4.7749714687606275}{1*0.85*cos(t)+0*0.85*sin(t)+2.34|0*0.85*cos(t)+1*0.85*sin(t)+1.34}
\parametricplot[linewidth=1.0pt,linestyle=dashed,dash=3pt 2pt]{-1.5125031955649737}{0.0989999773941644}{1*1.87*cos(t)+0*1.87*sin(t)+2.28|0*1.87*cos(t)+1*1.87*sin(t)+2.36}
\parametricplot[linewidth=1.0pt,linestyle=dashed,dash=3pt 2pt]{2.677755480453707}{4.9277297505381314}{1*0.54*cos(t)+0*0.54*sin(t)+1.39|0*0.54*cos(t)+1*0.54*sin(t)+3.02}
\parametricplot[linewidth=1.0pt,linestyle=dashed,dash=3pt 2pt]{1.8843071408090328}{2.583159724740975}{1*1.8*cos(t)+0*1.8*sin(t)+2.44|0*1.8*cos(t)+1*1.8*sin(t)+2.31}
\parametricplot[linewidth=1.0pt,linestyle=dashed,dash=3pt 2pt]{0.17396397227159288}{1.8107485743137344}{1*1.85*cos(t)+0*1.85*sin(t)+2.32|0*1.85*cos(t)+1*1.85*sin(t)+2.22}
\psdots[dotsize=6pt 0](1,3)
\psdots[dotsize=6pt 0](1,2)
\psdots[dotsize=6pt 0](1,1)
\psdots[dotsize=6pt 0](2,3)
\psdots[dotsize=6pt 0](2,2)
\psdots[dotsize=6pt 0](2,1)
\psdots[dotsize=6pt 0](3,3)
\psdots[dotsize=6pt 0](3,2)
\psdots[dotsize=6pt 0](3,1)
\psdots[dotsize=8pt 0,dotstyle=x](1.5,1.5)
\psdots[dotsize=8pt 0,dotstyle=x](2.5,1.5)
\psdots[dotsize=8pt 0,dotstyle=x](1.5,2.5)
\psdots[dotsize=8pt 0,dotstyle=x](2.5,2.5)
\psdots[dotsize=8pt 0,dotstyle=x](4.15,2.54)
\end{pspicture*}\hspace*{-2cm}\psset{xunit=1.0cm,yunit=1.0cm,algebraic=true,dotstyle=*,dotsize=3pt 0,linewidth=0.8pt,arrowsize=3pt 2,arrowinset=0.25}
\begin{pspicture*}(-0.36,-0.66)(6.08,3.82)
\psline[linewidth=1.6pt](2,3)(3,3)
\psline[linewidth=1.6pt](2,2)(2,3)
\psline[linewidth=1.6pt](1,3)(1,2)
\psline[linewidth=1.6pt](1,2)(2,2)
\psline[linewidth=1.6pt](2,1)(1,1)
\psline[linewidth=1.0pt,linestyle=dashed,dash=3pt 2pt](2.5,2.5)(2.5,1.5)
\psline[linewidth=1.0pt,linestyle=dashed,dash=3pt 2pt](2.5,1.5)(1.5,1.5)
\parametricplot[linewidth=1.0pt,linestyle=dashed,dash=3pt 2pt]{0.1852913062710569}{2.989280090214289}{1*1.34*cos(t)+0*1.34*sin(t)+2.83|0*1.34*cos(t)+1*1.34*sin(t)+2.3}
\parametricplot[linewidth=1.0pt,linestyle=dashed,dash=3pt 2pt]{1.1350083349985625}{2.0596018113318753}{1*1.84*cos(t)+0*1.84*sin(t)+3.37|0*1.84*cos(t)+1*1.84*sin(t)+0.87}
\parametricplot[linewidth=1.0pt,linestyle=dashed,dash=3pt 2pt]{4.6717617082515}{5.883637174045906}{1*1.71*cos(t)+0*1.71*sin(t)+2.57|0*1.71*cos(t)+1*1.71*sin(t)+3.21}
\parametricplot[linewidth=1.0pt,linestyle=dashed,dash=3pt 2pt]{1.5028552028831879}{2.100187227182034}{1*0.86*cos(t)+0*0.86*sin(t)+1.44|0*0.86*cos(t)+1*0.86*sin(t)+0.64}
\parametricplot[linewidth=1.0pt,linestyle=dashed,dash=3pt 2pt]{2.311045279311494}{4.0371577644321395}{1*0.54*cos(t)+0*0.54*sin(t)+1.37|0*0.54*cos(t)+1*0.54*sin(t)+0.99}
\parametricplot[linewidth=1.0pt,linestyle=dashed,dash=3pt 2pt]{-2.1993684311435224}{0.19067988888988038}{1*1.98*cos(t)+0*1.98*sin(t)+2.2|0*1.98*cos(t)+1*1.98*sin(t)+2.17}
\parametricplot[linewidth=1.0pt,linestyle=dashed,dash=3pt 2pt]{2.8453007051607155}{5.540388001593068}{1*0.52*cos(t)+0*0.52*sin(t)+3|0*0.52*cos(t)+1*0.52*sin(t)+1.35}
\parametricplot[linewidth=1.0pt,linestyle=dashed,dash=3pt 2pt]{5.428195851138397}{6.21778798060977}{1*2.24*cos(t)+0*2.24*sin(t)+1.91|0*2.24*cos(t)+1*2.24*sin(t)+2.69}
\psdots[dotsize=6pt 0](1,3)
\psdots[dotsize=6pt 0](1,2)
\psdots[dotsize=6pt 0](1,1)
\psdots[dotsize=6pt 0](2,3)
\psdots[dotsize=6pt 0](2,2)
\psdots[dotsize=6pt 0](2,1)
\psdots[dotsize=6pt 0](3,3)
\psdots[dotsize=6pt 0](3,2)
\psdots[dotsize=6pt 0](3,1)
\psdots[dotsize=8pt 0,dotstyle=x](1.5,1.5)
\psdots[dotsize=8pt 0,dotstyle=x](2.5,1.5)
\psdots[dotsize=8pt 0,dotstyle=x](1.5,2.5)
\psdots[dotsize=8pt 0,dotstyle=x](2.5,2.5)
\psdots[dotsize=8pt 0,dotstyle=x](4.15,2.54)
\end{pspicture*}}
\caption[Dual graph and dual RC state]{Left: The graph $G_3$ (solid) 
and its dual (dashed). Right: A RC state on $G_3$ (solid) and its 
dual configuration (dashed)}
\label{fig-dual}
\end{figure}

The \emph{random cluster model} (also known as the FK-model) was 
introduced by Fortuin and Kasteleyn in \cite{FK} to study lattice 
spin systems with a graph structure. It is defined on a graph 
$G=(V,E)$ by its state space 
$\O_{\rm RC}=\{\omega: \omega\subseteq E\}$ and the RC measure
\[
\mu_p(\omega) \;=\; 
\frac1Z\,p^{\abs{\omega}}\,(1-p)^{\abs{E}-\abs{\omega}}\,2^{C(\omega)},
\]
where $p\in(0,1)$, $Z$ is the normalization constant and 
$C(\omega)$ is the number of connected components in the graph 
$(V,\omega)$.
For a detailed introduction and related topics see the book 
\cite{G1}.\\
There is a tight connection between the Ising model and the 
random cluster model. Namely, if we set $p=1-e^{-\beta}$, 
we can translate an Ising configuration $\sigma\sim\pi_\beta$ to 
a random cluster state $\omega\sim\mu_p$ and vice versa.
To get an Ising configuration $\sigma\in\O_{\rm IS}$ from 
$\omega\in\O_{\rm RC}$ assign 
independent and uniformly random spins to each connected component 
of $\omega$. For the reverse way include all edges $e=\{e_1,e_2\}\in E$ 
with $\sigma(e_1)=\sigma(e_2)$ to $\omega$ with probability $p$. 
For details see \cite{ES}.\\
Therefore sampling an Ising configuration according to $\pi_\beta$ is 
equivalent to sampling a RC state from $\mu_p$ whenever both models are 
defined on the same graph $G$ and $p=1-e^{-\beta}$.\\
Another important concept in connection with the RC model is the 
duality of graphs (see e.g. \cite{G2}). 
Let $G=(V,E)$ be a finite, planar graph, i.e. 
without intersecting edges if we draw it in the plane (like our $G_L$). 
The \emph{dual graph} $G^*=(V^*,E^*)$ of $G$ is constructed as follows. 
Put a vertex in each face (including the infinite outer one) of the 
graph and connect 2 vertices by a edge if and only if the corresponding 
faces of $G$ share a boundary edge. It is clear, that the number of 
vertices can differ in the dual graph, but we have the same number of 
edges.\\
Additionally we define a \emph{dual configuration} 
$\omega^*\subseteq E^*$ in $G^*$ to a RC state $\omega\subseteq E$ in 
$G$ by
\[
e\in\omega \;\Longleftrightarrow\; e^*\notin\omega^*,
\]
where $e^*$ is the edge in $E^*$ that ``crosses'' $e$. (By the 
construction, this edge is unique.) 
See Figure \ref{fig-dual} 
for the graph $G_L$ with $L=3$ and its dual graph $G_L^*$ 
together with 2 corresponding RC states.

Now we can state the following theorem about the relation of the 
distribution of a RC state and its dual, see \cite{G2}.

\begin{proposition}{\cite[p.~164]{G2}}\label{prop-RC-dual}
Let $G=(V,E)$ be a finite, planar graph and $\mu_p$ be the random 
cluster measure on $G$. Furthermore let $G^*=(V^*,E^*)$ be the dual 
graph of $G$ and $\mu^*_{p^*}$ be the random cluster measure on $G^*$.\\
Then 
\[
\omega\sim\mu_p \;\;\Longleftrightarrow\;\; \omega^*\sim\mu^*_{p^*},
\]
where
\begin{equation}
p^* \;=\; 1 \,-\, \frac{p}{2-p}.
\label{eq-dual-p}
\end{equation}
\end{proposition}

Obviously, $(p^*)^*=p$. 
By Proposition \ref{prop-RC-dual} one can see that sampling 
from $\mu_p$ and sampling from $\mu^*_{p^*}$ is equivalent. 
It is straightforward to get the following Proposition.

\begin{proposition}\label{prop-Ising-dual}
Sampling from the Boltzmann distribution $\pi^G_\beta$ is 
equivalent to sampling from the Boltzmann distribution 
$\pi^{G^*}_{\beta^*}$, where
\begin{equation}\label{eq-dual-beta}
\beta^* \;=\; \log\left(\coth\,\frac{\beta}{2}\right).
\end{equation}
Additionally,
\[
\beta \,>\, \beta_c \;\;\Longleftrightarrow\;\; \beta^* \,<\, \beta_c.
\]
\end{proposition}

\begin{proof}
The equivalence was shown by the above procedure, i.e. if we want to 
sample from $\pi_\beta^G$, we can sample from $\pi^{G^*}_{\beta^*}$ 
generate a RC state with respect to $\mu^*_{p^*}$, go to the dual 
lattice with measure $\mu_p$ and finally generate a state 
according to $\pi^G_\beta$.
Since $p^{(*)}=1-e^{-\beta^{(*)}}$, the formula for $\beta^*$ comes 
from
\[
\beta^* \;=\; -\log(1-p^*) 
\;\overset{\eqref{eq-dual-p}}{=}\; \log\left(\frac{2-p}{p}\right) 
\;=\; \log\left(\coth\,\frac{\beta}{2}\right).
\]
This proves the statement.
\end{proof}


\section{Efficient sampling for the Ising model}\label{sec-efficient}

In this section we show an efficient algorithm to sample exactly 
from the Boltzmann distribution. 
But, before we prove that it is efficient, we state our 
sampling algorithm.\\
Therefore we first have to explain how the 
graph $G_L^*$ looks like. It is easy to obtain 
(see Figure \ref{fig-dual}) that 
$G_L^*=(V_L^*,E_L^*)$ is also a square lattice with $(L-1)^2$ 
vertices and an additional auxiliary vertex $v^*$, which is connected 
to every vertex on the boundary of it. We denote the operation of 
adding a vertex to a graph and connect it to all boundary vertices by 
$\cup_b$. So $G_L^*=G_{L-1}\cup_b v^*$.

\vspace*{1mm}
\begin{algorithm}[H]
\caption{\quad Sampling from the Ising model on the square lattice}
\begin{algorithmic}[1]
\Statex\Call{\bf Input:}{} An integer $L$ and the value of $\beta$
\Statex\Call{\bf Output:}{} An Ising configuration $\sigma\sim\pi^{G_L}_\beta$
\vspace{2mm}
\If{$\beta\le\beta_c$}
\vspace{1mm}
\State $\sigma=\text{CFTP}(G_L,\beta)$
\vspace{1mm}
\Else
\vspace{1mm}
\State $\widetilde\sigma=\text{CFTP}^+(G_{L-1},\beta^*)$, where 
				$\beta^*$ is given in \eqref{eq-dual-beta}
\State Define a Ising configuration $\sigma^*$ on 
				$G_L^*=G_{L-1}\cup_b v^*$ by 
\Statex\qquad\quad			$\sigma^*(v)=\widetilde\sigma(v)$ on 
							$V(G_{L-1})$ and $\sigma^*(v^*)=1$.
\State Generate a RC state $\omega^*$ from $\sigma^*$
\State Take the dual RC state $\omega=(\omega^*)^*$
\State Generate an Ising configuration $\sigma$ from $\omega$
\vspace{1mm}
\EndIf
\vspace{1mm}
\State \textbf{return} $\sigma$
\end{algorithmic}
\end{algorithm}

\begin{theorem}\label{th-main}
Let $G_L$ be the square lattice with $N=L^2$ vertices. Then, the 
algorithm from above outputs an exactly distributed Ising 
configuration with respect to $\pi_\beta^{G_L}$ in expected time 
smaller than
\begin{itemize}
	\item \quad $c_\beta\, N \,(\log N)^2$ 
		\qquad for $\beta\neq\beta_c=\log(1+\sqrt{2})$ and some $c_\beta>0$
	\vspace{2mm}
	\item \quad $16\, N^C \log N$ 
		\qquad
		for $\beta=\beta_c$, where $C$ is given in \eqref{eq-critical-C}.
\end{itemize}
\end{theorem}

\vspace{3mm}
\begin{proof}
The running time of the algorithm follows directly from Theorems 
\ref{th-mix-high} and \ref{th-mix-critical} and Prop.~\ref{prop-cftp}.
Therefore we only have to prove that the output $\sigma$ of the 
algorithm has the right distribution. In the case of 
$\beta\le\beta_c$ this is obvious. For $\beta>\beta_c$ we know from 
Proposition \ref{prop-Ising-dual} that $\sigma\sim\pi_\beta^{G_L}$, 
if the dual configuration $\sigma^*$ on $G_L^*$ 
(line 5 of Algorithm 2) is distributed according to 
$\pi_{\beta^*}:=\pi^{G_L^*}_{\beta^*}$. 
But by the construction of lines 4 and 5 of 
Algorithm 2, this is true. For this, note that 
$\pi_\beta(\eta)=\pi_\beta(-\eta)$ for all $\eta\in\Omega_{\rm IS}$. 
We get that for each vertex $v\in V$ (especially for $v^*$)
\[\begin{split}
\pi_\beta(\eta) \;&=\; \pi_\beta\bigl(\eta \cap \{\sigma\!:\sigma(v)=1\}\bigr) 
	\;+\; \pi_\beta\bigl(\eta \cap \{\sigma\!:\sigma(v)=-1\}\bigr) \\
&=\; \pi_\beta\bigl(\{\sigma\!:\sigma(v)=1\}\bigr)\;
			\pi_\beta\bigl(\eta \;\rule[-1.5mm]{0.2mm}{5mm}\; 
				\{\sigma\!:\sigma(v)=1\}\bigr) \\
&	\qquad\;+\; \pi_\beta\bigl(\{\sigma\!:\sigma(v)=-1\}\bigr)\;
			\pi_\beta\bigl(\eta \;\rule[-1.5mm]{0.2mm}{5mm}\; 
				\{\sigma\!:\sigma(v)=-1\}\bigr) \\
&=\; \frac12\,\Bigl[\pi_\beta\bigl(\eta \;\rule[-1.5mm]{0.2mm}{5mm}\; 
				\{\sigma\!:\sigma(v)=1\}\bigr) \;+\; 
				\pi_\beta\bigl(\eta \;\rule[-1.5mm]{0.2mm}{5mm}\; 
				\{\sigma\!:\sigma(v)=-1\}\bigr)\Bigr]\\
&=\; \frac12\,\pi_\beta\bigl(\{\eta,-\eta\} \;
				\rule[-1.5mm]{0.2mm}{5mm}\; \{\sigma\!:\sigma(v)=1\}\bigr).
\end{split}\]
The last equality comes from the fact that 
\[
\pi_\beta\bigl(\eta \;\rule[-1.5mm]{0.2mm}{5mm}\; 
				\{\sigma\!:\sigma(v)=-1\}\bigr) 
\;=\; \pi_\beta\bigl(-\eta \;\rule[-1.5mm]{0.2mm}{5mm}\; 
				\{\sigma\!:\sigma(v)=1\}\bigr). 
\]
Therefore we can sample from $\pi_\beta$ on $G_L^*$ by sampling $\eta$
from the conditional measure 
$\pi_\beta\bigl(\cdot \;\rule[-1.5mm]{0.2mm}{5mm}\;
\{\sigma\!:\sigma(v^*)=1\}\bigr)$ and then choose with probability $\frac12$
either $\eta$ or $-\eta$. 
If we now use that $G_L^*=G_{L-1}\cup_b v^*$ one can see that 
sampling on $G_L^*$ with respect to 
$\pi_\beta\bigl(\cdot \;\rule[-1.5mm]{0.2mm}{5mm}\;
\{\sigma\!:\sigma(v^*)=1\}\bigr)$ is the same as sampling $\widetilde\sigma$ from 
$\pi^{G_{L-1},+}_\beta$ and setting
\[
\sigma(v) \;=\; \begin{cases}
\widetilde\sigma(v), & v\in V(G_{L-1}) \\
1,			& v=v^*.
\end{cases}\]
Note that we omit the step of choosing $\sigma$ or 
$-\sigma$ with probability $\frac12$, because the RC state 
that will be generated would be the same.\\
This completes the proof.\\
\end{proof}

\begin{remark}
Note that the same technique works also for the $q$-state Potts model.
This model consists of the state space $\O_{\rm P}=\{1,\dots,q\}^V$ 
and the same measure $\pi_\beta$. 
In this case we consider the random cluster measure
\[
\mu_{p,q}(\omega) \;=\; 
\frac1Z\,p^{\abs{\omega}}\,(1-p)^{\abs{E}-\abs{\omega}}\,q^{C(\omega)}
\]
and the connection of the models is again given by $p=1-e^{-\beta}$.\\
A recent result of Beffara and Duminil-Copin \cite{BDC} shows that 
the self-dual point of the RC model corresponds to the critical 
temperature of the Potts model $\beta_c(q)=\ln(1+\sqrt{q})$ in the 
same way as in the case $q=2$ (i.e. the Ising case). 
Therefore, a sampling algorithm for the 
Potts model above (and at) the critical temperature is enough to 
sample at all temperatures.
\end{remark}

\linespread{1}
\bibliographystyle{amsalpha}
\bibliography{Bibliography}


\end{document}